\documentclass[11pt]{article}

\usepackage{amssymb,amsmath,amsfonts,amsthm,comment}
\usepackage{latexsym}
\usepackage{graphics}
\usepackage{indentfirst}

\newtheorem{theorem}{Theorem}
\newtheorem{lemma}{Lemma}

\newtheorem{prob}[theorem]{Problem}

\def\qed
 {\ifhmode\unskip\nobreak\hfill$\Box$\medskip\fi
 \ifmmode\eqno{\Box}\fi}

\normalsize

\title{{\bf \huge  List colorings with distinct list sizes, \\the case of complete bipartite graphs}}
\author{{\bf Zolt\'an F\"uredi}${}^{1}$\\ %[-0.8ex]
\small R\'enyi Institute of the Hungarian Academy\\[-0.8ex]
\small Budapest, P.O.Box 127, Hungary, H-1364,\\[-0.8ex]
\small \texttt{furedi@renyi.hu}\\[-0.8ex]
\small and\\[-0.8ex]
\small Department of Mathematics\\[-0.8ex]
\small University of Illinois at Urbana-Champaign\\[-0.8ex]
\small Urbana, IL 61801, USA \\[-0.8ex]
\small \texttt{z-furedi@illinois.edu}
\and {\bf Ida Kantor}\\
\small Institute for Theoretical Computer Science\\[-0.8ex]
\small Charles University\\[-0.8ex]
\small Malostransk\'{e} n\'{a}m. 25, Praha 1, Czech Republic \\[-0.8ex]
\small \texttt{ida@kam.mff.cuni.cz}}
\date{${}$}

\begin{document}

\maketitle

\renewcommand{\thefootnote}{\empty}
\footnotetext{
%05D15  Transversal (matching) theory
%05B40  Packing and covering
\hfill    {\tt  [{\jobname}.tex]}\newline
Submitted to  \emph{???}   \hfill
Printed on \today \newline %\indent\indent
${}^1$
Research supported in part by the Hungarian National Science Foundation
OTKA, by the National Science Foundation grant DMS 09-01276,
and by a European Research Council Advanced Investigators Grant 267195.\\
Acknowledgement:
This research was supported in part by the Institute for Mathematics and its Applications
with funds provided by the National Science Foundation.\\
${}$\quad
A preliminary version was presented {\em in:}
% {\sc Z. F\"uredi and I. Kantor:}
% List colorings with distinct list sizes, the case of complete bipartite graphs,
 European Conference on Combinatorics, Graph Theory and Applications (EuroComb 2009), 323--327,
Electron. Notes Discrete Math. 34, Elsevier Sci. B. V., Amsterdam, 2009.
}

\vskip -1cm

\begin{abstract}
Let $f:V  \rightarrow \mathbb{N}$ be a  function on the vertex set of the graph $G=(V,E)$.
The graph $G$ is {\em $f$-choosable} if for every collection of lists with list sizes specified by $f$
 there is a proper coloring using colors from the lists.
The sum choice number, $\chi_{sc}(G)$, is the minimum of $\sum f(v)$, over all functions $f$ such that $G$ is $f$-choosable.
It is known (Alon 1993, 2000) that if $G$ has average degree $d$, then
 the usual choice number $\chi_\ell(G)$ is at least $\Omega(\log d)$, so they grow simultaneously.

In this paper we show that $\chi_{sc}(G)/|V(G)|$ can be bounded while the minimum degree
$\delta_{\min}(G)\rightarrow \infty$.
%(This is not true for the list chromatic number, $\chi_\ell(G)$).
Our main tool is to give tight estimates for the sum choice number of %for
 the unbalanced complete bipartite graph $K_{a,q}$.
\end{abstract}

\titlepage

\section{Average list sizes and planar graphs}

Given a graph $G$ and a list of colors $L(v)$ for each vertex $v\in V(G)$, we say that
  $G$ is $L$-{\em choosable} (or that $L$ is {\em sufficient}) if it is possible to choose $c(v)\in L(v)$
  for all $v$ so that $c:V(G)\to \cup L(v)$ is a proper coloring of $G$.
The {\em choice number} (or {\em list chromatic number}) $\chi_\ell$ is the minimum $t$
  such that for every assignment $L$ with $|L(v)|\geq t$ for all $v\in V$,
  the graph is $L$-choosable.
It is well-known (Thomassen~\cite{thom:ever94}) that
\begin{align}\label{eq11}
 \chi_\ell(P)\leq 5
\end{align}
for every planar graph $P$, and this is the best possible (\cite{voig:list93}).

However, if we allow distinct list sizes, then the average size can be smaller.
For example, Thomassen's beautiful proof for~(\ref{eq11}) gives that
 if $P$ is an $n$-vertex planar graph, $v_1, \dots, v_t$ are its external vertices (in this order)
  and the list sizes are
\begin{align}\label{eq33}
|L(v)|=\left\{ \begin{array}{ll}
         1 & \mbox{ for $v=v_1$},\\
         2 & \mbox{ for $v=v_2$},\\
         3 & \mbox{ for $v=v_3, \dots, v_t$},\\
         5 & \mbox{ for the inner vertices},
         \end{array} \right. \end{align}
then $P$ is $L$-choosable.

Consider a function $f:V(G)\rightarrow \mathbb{N}$.
An {\em $f$-assignment} is an assignment of lists $L(v)$ to the vertices
$v\in V(G)$ such that $|L(v)|=f(v)$ for all $v$.
The function $f$ is {\em sufficient} %c
if $G$ is $L$-choosable for all $f$-assignments $L$.
We define the {\em sum choice number} of $G$, denoted by $\chi_{sc}(G)$,
 as the minimum of $\sum_{v\in V(G)}f(v)$ over all sufficient %c 
$f$.

Sum choice numbers were introduced by Isaak~in~\cite{isaa:suml02}
 who proved that if $G$ is the line-graph of $K_{2,q}$ then
$\chi_{\rm sc}(G)=q^2 + \lceil 5q/3\rceil$.
Various classes of graphs were investigated by
 Isaak~in~\cite{isaa:suml04}, by
 Berliner, Bostelmann, Brualdi, and Deaett~\cite{berl:suml06} and by Heinold in~\cite{hein:suml06} and~\cite{hein:sumc09}.

Thomassen's theorem (\ref{eq33}) implies that $\chi_{sc}(P)\leq 5n-9$ for planer $P$ ($n\geq 2$).
In fact, more is true. It is easy to show   (see, e.g.,~\cite{isaa:suml04}) that for every graph
\begin{align}\label{eq44}\chi_{sc}(G)\leq |V(G)|+|E(G)|
\end{align}
 holds.
Hence $\chi_{sc}(P)\leq 4n-6$.
Our first result is a slight improvement.

\begin{theorem}\label{t11}\enskip
 Let $P$ be an $n$-vertex planar graph.
 There exists an $f:V(P)\rightarrow \mathbb{N}$ such that
% {\rm (1)}
$\sum f(v)\leq 4n-6$,
% {\rm (2)}
$\max f(v)\leq 6$, and
% {\rm (3)}
$P$ is $f$-choosable. \qed
\end{theorem}

{\em Proof.} Consider a linear order of the vertices of $P$, and let $\hat{d}(v)$ be the number of neighbors of $v$ that precede it.
The function $f(v)=\hat{d}(v)+1$ is a sufficient %c
 function, so $\sum (\hat d(v) +1)$
% every ordering of $V(P)$
  yields an upper bound on $\chi_{sc}(P)$.
Since every planar graph has a vertex of degree at most 5, it is possible to order the vertices so that $\hat{d}(v)\leq 5$ for all $v$.
\hfill{$\Box$}

\section{Unbalanced complete bipartite graphs}

Erd\H os, Rubin and Taylor (see, e.g.,~\cite{alon:rest93}) showed for the complete
 bipartite graph that
\begin{align}\label{eq22}
 \chi_\ell(K_{q,q})=\Theta (\log q).
\end{align}

If one of the parts is substantially smaller than the other one,
 then allowing different list sizes %again
  results in smaller average lists.
It is easy to show $\chi_{sc}(K_{1,q})=2q+1$ (as for every tree on $q+1$ vertices).
Berliner, Bostelmann, Brualdi, and Deaettet~\cite{berl:suml06} showed that
for all $q\geq 1$ we have
\begin{equation}\label{t1}
  \chi_{sc}(K_{2,q})=2q+1+\lfloor \sqrt{4q+1}\rfloor,
  \end{equation}
and Heinold~\cite{hein:sumc09} proved
\begin{equation}\label{t2}
 \chi_{sc}(K_{3,q})=2q+1+\lfloor \sqrt{12q+4}\rfloor.
 \end{equation}
Our main result % extends Theorems~\ref{t1} and~\ref{t2} to
  deals with the sum choice number of $K_{a,q}$ with arbitrary $a$.

\begin{theorem}\label{main}
 There exist positive constants $c_1$ and $c_2$ such that for all $a\geq 2$ and $q\geq 4 a^2\log  a$
\[ 2q+c_1a\sqrt{q\log {a}}\leq \chi_{sc}(K_{a,q})\leq 2q+c_2a\sqrt{q\log {a}}. \]
\end{theorem}

It is known 
 that $\chi_\ell$ is not independent of the average degree.
Alon~\cite{alon:rest93,A2000} proved that
 for some constant $c>0$, every graph $G$ with average degree $d$ has
\begin{equation}\label{alon}
 \chi_\ell({G})\geq c \log d.
\end{equation}

An easy %One of the most interesting
 corollary of our Theorem~\ref{main} is %(\ref{eq55}),
 that if different list sizes are allowed, then such dependence does not exist.
Indeed,
we have
\begin{align*}%\label{eq55}
   \lim_{a\to \infty, \atop q>> a^2 \log a} \frac{2|E(K_{a,q})|}{{a+q}}=\infty, \quad\quad
   \lim_{a\to \infty, \atop q>> a^2 \log a} \frac{\chi_{sc}({K_{a,q}})}{{a+q}}=2.
\end{align*}
So the structure of the graph plays a more important role in determining
 the sum choice number than in the case of the list chromatic number.

\section{Upper bound, there are sufficient %c
 short lists}

Throughout this paper, the two parts of the complete bipartite graph
 $K_{a,q}$ will be denoted by $A$ and $Q$, with $|A|=a$ and $|Q|=q$.
\begin{theorem}
\label{ub}
Suppose that $a,q\in \mathbb{N}$ with $q \geq a \geq 2$. Then
\begin{align*}
\chi_{sc}(K_{a,q})\leq 2q+ a\lceil\sqrt{32 q (1+\log a)}\rceil.
\end{align*}
\end{theorem}

{\em Proof.}
To prove the upper bound, we present a function $f$ with
   $\sum_{v\in A\cup Q}f(v) = 2q+a\lceil\sqrt{32 q (1+\log a)}\rceil$
   %%% !!! = AND NOT \geq
    such that every $f$-assignment is sufficient. %c

Define $f$ as
\[f(v)=\left\{ \begin{array}{ll}
         r & \mbox{for $v\in A$};\\
         2 & \mbox{for $v\in Q$}\end{array} \right. \]
where $r$ will be defined later in (\ref{eq3}) as any integer
  $r\geq \sqrt{32 q (1+\log a)}$.
Let $L$ be an arbitrary $f$-assignment, i.e., $|L(v)|=f(v)$ for all $v$.

\def\LA{{{\mathcal  L}_A}}
\def\LQ{{{\mathcal L}_Q}}

Consider $S:=\bigcup_{v\in A\cup Q} L(v)$.
The assignment $L$ yields a (multi)hypergraph and a multigraph on the same
 vertex set $S$ and with edge sets  $\LA:=\{ L(u): u\in A\}$ and
 $\LQ:=\{ L(v): v\in Q\}$, respectively.
Sufficiency %c
 of $L$ means that one can find a set $T\subset S$
 meeting all hyperedges of $\LA$ such that $S\setminus T$ meets all edges of $\LQ$,
so $T$ is an independent set in the graph $\LQ$.
Given $T$ the choice function $c$ can be defined as
$$c(u)\in L(u)\cap T,\text{ for }u\in A$$
and
$$c(v)\in L(v)\cap (S\setminus T),\text{ for }v\in Q.
 $$
We are going to construct such $T$  by a 2-step random process.

Let us pick, randomly and independently, each element of $S$
 with probability $p$.
Let $B$ be the random set of all elements picked.
Define a random variable $X_u$ for each  $u\in A$ as $X_u=|L(u)\cap B|$,
 and the random variable $Y$ by
 $$ Y:=| \{ v\in Q: L(v)\subseteq B \}|, $$
 so $Y$ is the number of edges of $\LQ$ spanned by $B$.
Remove an element $c(v)\in L(v)$ for each edge of $\LQ$ spanned by $B$,
 the remaining set $T\subset B$ is certainly independent in $\LQ$, and
if $Y< X_u$ for each $u\in A$, then $T$ meets all $L(u)\in \LA$ and
we are done.

\def\Pr{{\rm Prob \,}}
The expected value of $Y$ is $p^2q$, so  Markov inequality gives
\begin{align}\label{eq1}
 \Pr(Y< 2p^2q)\geq \frac{1}{2}.
\end{align}

The expected value of $X_u$ is $pr$, so  Chernoff inequality gives
\begin{align*}
 \Pr(X_u< EX_u-t) < e^{-t^2/2rp},
\end{align*}
for any $t>0$.
Hence
\begin{align}\label{eq2}
 \Pr(X_u\geq pr-t \mbox{ for all $u\in A$}) > 1- ae^{-t^2/2rp},
\end{align}
and this is at least $1/2$ for $t^2\geq 2rp \ln(2a)$.
The sum of probabilities in (\ref{eq1}) and (\ref{eq2}) is larger than 1, so there is
 an appropriate choice of $B$ (and then $T$) if
 $t^2=2rp(1+\log a)$ and $pr-t \geq 2p^2q$.
We can define, e.g.,
\begin{align}\label{eq3}
 p:= \sqrt{\frac{2(1+\log a)}{q}} \quad{\rm and}\quad r\geq 4pq=\sqrt{32(1+\log a)q}.
  \end{align} \hfill{$\Box$}

%\newpage

%%%%%%%%%%%%%%%%%%%%%%%%%%%%%%%%%%%%%%%%%%%%%%%%%%%%%%%%%%%%%%%%%%%%%%%%%%%%%%%%%%%%%%%%%%%%%%%%%%%%%55
\section{Lower bound, much shorter lists are not sufficient %c
}\label{lower}

To prove that $\chi_{sc}(G)\geq k$ for a particular $k$,
 we need to show that for every $f$ with $\sum_{v\in G}f(v)=k$, there exists an insufficient %c
 $f$-assignment.
First, we show how to construct an insufficient %c
 assignment for some special $f$.

\begin{lemma}\label{constr}
Let $t\geq 2$ and $\ell\geq 1$. For $a=2^t$ and $q=t\ell^2$, there exists an insufficient %c
 assignment $L$ with
\[ |L(v)|=\left\{ \begin{array}{ll}
        t\ell & \mbox{for $v\in A$}, \\
        2 & \mbox{for $v\in Q$}.\end{array} \right. \]
\end{lemma}

{\em Proof.}
Take $2t$ pairwise disjoint sets $X_i,Y_i$ of size $\ell$, $Z=\cup (X_i\cup Y_i)$.
Identify the elements of $A$ by the set of $0$-$1$ vectors of length $t$,
 $A= \{ 0,1 \}^t$.
For a vector $(\varepsilon_1, \dots, \varepsilon_t)\in A$ define $L(v)$ as
 $(\bigcup_{\varepsilon_i=1} X_i) \cup (\bigcup_{\varepsilon_j=0}  Y_j )$.
So $L(v)$ contains either $X_i$ or $ Y_i$ for all $i$.
% Replace the vertices $u_i,v_i$ with disjoint vertex .
Let the graph $G$ be the union of $t$ complete bipartite graphs on the vertex set $Z$ by setting $E(G)=\bigcup_{i=1}^t (X_i\times Y_i)$
 and define the lists $L(v)$ for $v\in Q$ as the edges of $G$.
The number of edges of $G$ is $t\ell^2=q$, so a one-to-one mapping can be done.

Every independent set $T$ of  $G$ contains at most one vertex from each $X_i\cup Y_i$ so it
  cannot meet all hyperedges of ${\mathcal L}_A$, where ${\mathcal L}_A:= \{ L(v): v\in A\}$.
This means that this assignment $L$ is not sufficient. %c
\hfill{$\Box$}

Note that with this choice of $a$ and $q$, we have $|L(v)|= \sqrt{q\log_2 a}$ for $v\in A$.
Also notice that if we remove some elements from the lists in the above construction, the resulting list assignment is still insufficient. %c

\begin{theorem}\label{lb}
 If $a\geq 2$ and $q>4 a^2\log a$, then
\begin{align*}
 \chi_{sc}(K_{a,q})\geq 2q+ 0.06 a\sqrt{q\log a}
\end{align*}
\end{theorem}

\noindent{\em Proof.}
Suppose that $f:V(K_{a,q})\to \mathbb N$ with $\sum_{v\in A\cup Q} f(v)=2q+as$ where $s\leq 0.06 \sqrt{q\log a}$.
We will find an insufficient %c
 $f$-assignment. % Suppose for contradiction that no such assignment exists.

Let $q_1$, $q_2$ and $q_3$ be the numbers of vertices $v\in Q$ with $f(v)=1$, $f(v)=2$ and $f(v)\geq 3$, respectively.
If $f(u)\leq q_1$ for some $u\in A$, then $f$ is obviously insufficient. %c
From now on, we suppose that
  $f(u) > q_1$ for each $u\in A$. % (otherwise we can find a non-choosable assignment)
We obtain
\begin{align}\label{eq:q1+q3}
2q+as=\sum_{v\in A\cup Q}f(v)\geq aq_1+(q_1 +2q_2 +3q_3)\geq 2q +q_1+q_3.
\end{align}
It follows that $q_1+q_3\leq as$ and $Q$ has at least $q-as$ vertices with lists of size 2.
Let $q^*=q_2$ and let $a^*$ be the largest power of $2$ not exceeding $\frac{1}{2}a$.

If there are at least $a^*$ vertices $u\in A$ with $f(u)\leq \sqrt{q^*\log_2 a^*}$,
 then we can use Lemma~\ref{constr} %%% together with the subsequent remark
 to construct an insufficient %c
 assignment.

If this does not hold, %this is not true,
 then $A$ has more than $\frac{a}{2}$ vertices with lists of size greater than $\sqrt{{(q-as)\log_2 a^*}}$.
Using $a^*> a/4$ we obtain
\begin{align}\label{tooMany}
 \sum_{v\in A\cup Q}f(v)&\geq \frac{1}{2}a  \sqrt{{(q-as)\log_2 (a/4)}}+2q-as. % >2q+as.
\end{align}
The rest is a little calculation to show that here the right hand side exceeds
 $2q+ 0.06 a\sqrt{q\log a}$
 for $a\geq 5$, $q>4 a^2\log a$ and $s<  0.06 \sqrt{q\log a}$.
Finally, the case $a\leq 4$ (in fact $a\leq 30$) follow from (\ref{t1}) and (\ref{t2}), completing the proof.
\hfill{$\Box$}

% \vspace{3mm}
\bigskip
Let us remark that if we choose the constants in the proof of Theorem~\ref{ub} more carefully, we can improve the constant $\sqrt{32}$ to $3.67$.
Using a randomized construction, it is possible to improve the constant 0.06 in Theorem~\ref{lb} to % %%%%%%%%%%%%%%%%%%%%%%%%%%%%%%%%%%%%%%%%%%%%%%%%%%%%
\section{For fixed $a$, a limit exists as $q\rightarrow \infty$}

In this section we suppose that $a\geq 2$ is a fixed integer.
We have proved bounds for $\alpha_q:=(\chi_{sc}(K_{a,q})-2q)/{\sqrt{q}}$. Now we show that in fact the limit exists when $q$ tends to $\infty$.

\begin{theorem}\label{limExists}
For fixed $a$, the limit $\lim_{q\rightarrow \infty} \frac{\chi_{sc}(K_{a,q})-2q}{\sqrt{q}}$ exists.
\end{theorem}

First, we consider a simpler problem and consider only {\em type II} assignments of $K_{a,q}$ which means
 $f(v)=2$ for all $v\in Q$.
Define $\chi_{sc2}(K_{a,q})$ to be the minimum of $\sum_{A\cup Q}f(v)$
 where $f$ runs over all sufficient %c
 type II functions.
Obviously $\chi_{sc2}(K_{a,q})\geq \chi_{sc}(K_{a,q})$.

\begin{theorem}\label{lE2}
For fixed $a$, the limit $\lim_{q\rightarrow \infty} \frac{\chi_{sc2}(K_{a,q})-2q}{\sqrt{q}}$ exists.
\end{theorem}

A type II $f$ is not sufficient %c
 if and only if there exists a % non-choosable assignment $(\mathcal{L},G)$ for $f$ and $q$, i.e.
 a hypergraph $\mathcal{L}$ with edges $L_i$ satisfying $|L_i|=f_i$ for $i=1,\dots,a$ and a graph $G$
 on $V(\mathcal{L})$ with at most $q$ edges, such that no transversal of $\mathcal{L}$ is an independent set in $G$.

For $I\subseteq [a]$, define $X_I=\cap_{i\in I}L_i$. An insufficient %c
 II assignment is {\em symmetric} if for all pairs $I\neq J$,
 the bipartite subgraph of $G$ induced by $X_I$ and $X_J$ is either empty of complete, and for each $I$, $X_I$ induces the empty graph.
Without loss of generality we may assume that an insufficient %c
 type II assignment is symmetric, as the following lemma demonstrates.
From now on, in this section, all assignments are of type II, except when stated otherwise.

\begin{lemma}
Given $a$, $q$ and $f$ an insufficient %c
 type II assignment exists if and only if a symmetric insufficient %c
 assignment exists.
\end{lemma}

{\em Proof.} Suppose that $(\mathcal{L}, G)$ is an insufficient %c
 assignment.
If $u$ and $v$ belong to the same $X_I$, then no minimal transversal of $\mathcal{L}$ contains both of them.
We can therefore delete all edges induced by $X_I$.

Now suppose that $u,v\in X_I$ and $|N(u)|\leq |N(v)|$. Replace the neighborhood of $v$ by the neighborhood of $u$.
It is still true that every transversal of $\mathcal{L}$ induces an edge of $G$.
Repeated application of this procedure eventually produces a symmetric insufficient %c
 assignment. \hfill{$\Box$}

%%%\vspace{3mm}
\bigskip
{\em Proof of Theorem~\ref{lE2}.}
Consider a symmetric insufficient %c
 assignment $f$ for $K_{a,q}$.
Let $\mathcal{L},G$ and $X_I$ (for $I\subseteq [a]$) be as before, $x_I=|X_I|$.
%%% Replace each $x_I$ with a vertex $v_I$.
Let $V:=\{ v_I: I\subseteq [a]\}$ be a $2^a$-element set.
Let $R$ be the {\em reduced graph} of the symmetric insufficient %c
 assignment, i.e.,
 the graph with $V(R)=\{v_I;x_I\neq 0\}$ and whose edges correspond to the complete bipartite subgraphs of $G$.
% need only {\em nonempty} x_I to be the vertices of the graph and hypergraph!
Similarly, the hypergraph $\mathcal{L}$ turns into the reduced hypergraph on the same vertex set, $V(R)$.
The graph $R$ is {\em blocking}, %c  %find better word
i.e., every vertex cover of the reduced hypergraph contains %%%induces
 an edge of $R$.
The vector $x=(x_I)$ satisfies $\sum_{IJ\in E(R)}x_I x_J\leq q$, and $x_I=0$ whenever $v_I\not\in V(R)$.
The set $A_R^q$ of all such $x$ lies in the non-negative orthant of $\mathbb{R}^{2^{[a]}}$ and is bounded by a quadric surface which depends on $R$ and $q$.

Define the linear map $\varphi:\mathbb{R}^{2^{[a]}}\rightarrow \mathbb{R}^a$ by $\varphi(x)=(f_1,\dots,f_a)$ where $f_i=\sum_{i\in I}x_I$.
%Let $B'_{R,q}:=\varphi(A'_{R,q})$.
The function $f$ is insufficient %c
 for this $q$ if and only if $f$ is the image of some integer point $x$
%$f\in \cup A'_{R,q}$ where the union is taken over all nonchoosable $R$.
that is in $A_R^q$ for some blocking %c
 $R$.

If there exists an insufficient %c
 $f$-assignment for every integer vector $f$ such that $\sum f_i=k$, then we have $\chi_{sc2}(K_{a,q})-2q>k$.
We are therefore looking for the maximum $k$ such that every integer point on the hyperplane $\sum f_i=k$ is the image
 (under $\varphi$) of some integer point in $\bigcup A_R^q$, where the union is taken over all (but finitely many) blocking %c
 $R$'s.

Let us normalize everything by $\sqrt{q}$. For every blocking %c
$R$, define
$$A_R:=\{x: \sum_{IJ\in E(R)}x_I x_J\leq 1 \mbox{, and }x_I=0\mbox{ for }v_I\not\in V(R)\} \quad\mbox{ and } \quad B_R:=\varphi(A_R).$$ %remove B_R?
For every $R$ we now have only one quadric surface, independent of $q$.
We say that a vector $v$ is a {\em $q$-grid point} if $\sqrt{q}\cdot v$ is an integer point.

For every $q$, define $k_q$ to be the maximum $k$ such that every $q$-grid point on the hyperplane $\sum f_i=k$ is the image of some $q$-grid point in $\bigcup A_R$.
Also, define $\beta$ to be the maximum $k$ such that the simplex $C_k:=\{f: \sum f_i\leq k\}$ is a subset of $\bigcup B_R$.

We want to prove that the limit $\lim k_q$ exists and equals $\beta$. That is, we want to prove that for every $\varepsilon$, if $q$ is large enough,
\begin{itemize}
\item every $q$-grid point in $C_{\beta-\varepsilon}$ is the image under $\varphi$ of some $q$-grid point in $\bigcup A_R$, and
\item there is a $q$-grid point in $C_{\beta+\varepsilon}$ which is not the image of any $q$-grid point in $\bigcup A_R$.
\end{itemize}

% \vspace{3mm}
\bigskip
To prove the first claim, fix $q$ and let $f$ be a point on the hyperplane $\sum f_i=\beta$.
The point $f$ is in $\bigcup B_R$, so it is the image of some $x\in \bigcup A_R$.
Each set $A_R$ is a {\em downset} in the sense that with every $x$ it also contains all points $z$ such that $z_i\leq x_i$ for all $i$.
It follows that $y:=\frac{\lfloor \sqrt{q} \cdot x \rfloor}{\sqrt{q}}$ is a $q$-grid point in $\bigcup A_R$.
Each entry of $y$ differs by at most $\frac{1}{\sqrt{q}}$ from the corresponding entry of $x$,
 and a simple computation suffices to show that the distance of $\varphi(y)$ and $f$ is at most $\frac{c}{\sqrt{q}}$, where $c$ is a constant dependent only on $a$.
That is, for each point $f$ on the hyperplane $\sum f_i=\beta$ we have found, in distance at most $\varepsilon_q:= \frac{c}{\sqrt{q}}$,
 an image of a $q$-grid point from $\bigcup A_R$. Call this point $f'$.

Note that, by definition of $\chi_{sc}$, whenever a $q$-grid point $f$ is the image of a $q$-grid point in $\bigcup A_R$,
 the same is true for all $q$-grid points in the box $D_f:=\{g: g_i\leq f_i\mbox{ for all }i\}$.

Let $h$ be a $q$-grid point such that $\sum h_i\leq \beta-\varepsilon_q\cdot\sqrt{a}$.
Let $f$ be its perpendicular projection on the hyperplane $\sum f_i=\beta$ and find the corresponding $f'$.
Since the distance of $f$ and $f'$ is at most $\varepsilon_q$, the point $h$ belongs to $D_{f'}$, and hence it is the image of a $q$-grid point in $\bigcup A_R$.
Choosing $q$ large enough so that $\varepsilon_q\cdot\sqrt{a}\leq \varepsilon$ for our given $\varepsilon$ concludes the proof.

% \vspace{3mm}
\bigskip
Now we prove the second claim. Let $f$ be a point outside $\cup B_R$, but within the distance $\varepsilon$ from $C_\beta$.
Take a bounded cube $Q\subseteq \mathbb{R}^a$ that contains $f$.
Now take a bounded cube $S$ in $\mathbb{R}^{2^a}$ which contains all points $x$ such that $\varphi(x)\in Q$.
Then $T:=S\cap (\bigcup A_R)$ is a compact set, so $\varphi$ maps it to a compact set.
The complement of $\varphi(T)$ in $\varphi(S)$ is open, and contains $f$. Note that $(\bigcup B_R)\cap Q\subseteq \varphi(T)$.

Therefore, for some small $\delta$, the $\delta$-ball around $f$ is outside $\bigcup B_R$.
If $q$ is large enough, the ball contains some $q$-grid point.
This point not only has no $q$-grid preimages in $\cup A_R$, it has no preimages in $\bigcup A_R$ whatsoever, and the claim is proven. \hfill{$\Box$}

% \vspace{3mm}
\bigskip
{\em Proof of Theorem~\ref{limExists}.} Define sequences $\{\alpha_q\}_{q=1}^{\infty}$ and $\{\beta_q\}_{q=1}^{\infty}$ as follows
\[ \alpha_q=\frac{\chi_{sc} (a,q)-2q}{\sqrt{q}} \quad\mbox{ and }\quad \beta_q=\frac{\chi_{sc2} (a,q)-2q}{\sqrt{q}}.\]
It was already mentioned that $\alpha_q\leq \beta_q$ for all $q$.

The argument in the proof of Theorem~\ref{lb} shows that whenever we have an insufficient %c
 function $f$ for $K_{a,q}$,
 we can delete at most $d(q)=O(\sqrt q)$ vertices of $Q$ and get an insufficient %c
 function for $K_{a,q-d(q)}$ where $f(v)=2$ for $v\in Q$.
 We therefore have $\alpha_q\sqrt{q}=\chi_{sc}(a,q)-2q\geq \chi_{sc2}(q,q-d(q))-2(q-d(q))=\beta_{q-d(q)}\sqrt{q-d(q)}$.
 We get the following relationship between $\alpha_q$ and $\beta_q$:
\[\beta_q\geq \alpha_q\geq \frac{\sqrt{q-d(q)}}{\sqrt{q}} \beta_{q-d(q)}.\]
The limit $\lim_{q\rightarrow \infty}\beta_q$ exits by Theorem~\ref{lE2}. Since $d(q)=O(\sqrt{q})$,
we have \[\lim_{q\rightarrow \infty} \frac{\sqrt{q-d(q)}}{\sqrt{q}} \beta_{q-d(q)}=\lim_{q\rightarrow \infty}\beta_q,\] which proves the claim. \hfill{$\Box$}

%%%%%%%%%%%%%%%%%%%%%%%%%%%%%%%%%%%%%%%%%%%%%%%%%%%%5555555
\section{Graphs with large independent sets and a generalization of Tur\'an's theorem}

Let $G_{a,q}$ be the graph that we get from $K_{a,q}$ by inserting an edge $\{u,v\}$ for every pair of distinct $u,v\in A$.

\begin{theorem}\label{split} There exist positive constants $c_1$ and $c_2$, independent of $q$ and $a$, such that
 for $q> a\geq 2$ we have
\[2q+ c_1 a\sqrt{q(a-1)} \leq\chi_{sc}(G_{a,q})\leq 2q+ c_2 a \sqrt{q(a-1)}. \]
%2q+a\sqrt{2(a-1)q}\cdot(1+o(1))\]

\end{theorem}

We will use the following generalization of Tur\'{a}n theorem.
For given positive integers $s$ and $k$, let $t(s,k)=\min \sum_{1\leq i\leq k} {{d_i}\choose{2}}$, where the minimum is taken over all
 non-negative integer sequences $(x_1,\dots,x_k)$ such that $\sum d_i=s$.

\begin{theorem}\label{turan}
Let $s,a\geq 2$ be integers and let $G$ be a graph with  less than $t(s,a-1)$ edges.
Let $L_1,\dots,L_a\subset V(G)$ be sets of size $s$.
Then there exists a system of distinct representatives $\{u_1,\dots,u_a\}$ of $L_1,\dots,L_a$,
 (i.e., $u_i\in L_i$) which is an independent $a$-element set in $G$.
\end{theorem}

The case $V(G)=L_1=\dots=L_a$ gives (the dual form of)  Tur\'{a}n's theorem.

Let us also remark that the result concerning $|E(G)|$ is sharp: %if $G$ is allowed to have $t(s,a-1)$ edges,
 taking $L_1=\dots=L_a=V(G)$, $|V(G)|=s$  with $G$ being the disjoint union of $a-1$ cliques of almost equal sizes
 provide a graph of $t(s,a-1)$ edges and a family without any system of distinct representatives
 which is independent in $G$ (because $G$ has no any independent set of size $a$).

\bigskip
%%%\vspace{3mm}
{\em Proof.}
We define the sequence of distinct vertices $u_1,\dots,u_a$ one by one by an algorithm,
 such that $u_k\in L_{i_k}$, where $\{i_1,i_2,\dots,i_a\}$ is a permutation of $\{1,2,\dots,a\}$
  and also the set $\{u_1,\dots,u_a\}$ is independent in $G$.

Let $V_1=L_1\cup \dots \cup L_a$ and $G_1=G[V_1]$ (the restriction of $G$ to $V_1$).
Let $u_1$ be the vertex of minimum degree in $G_1$, $D_1$ the closed neighborhood of $u_1$ in $G_1$, and $d_1=|D_1|$.
Let $L_{i_1}$ be one of the hyperedges containing $u_1$.

If $u_j,D_j$ and $L_{i_j}$ are already defined for $j=1,\dots,k$, consider
$V_{k+1}=(\bigcup_{i\not\in \{i_1,\dots,i_k\}} L_i)\setminus (D_1\cup \dots \cup D_k)$, $G_{k+1}= G[V_{k+1}]$,
and let $u_{k+1}$ be a vertex of minimum degree in $G_{k+1}$, $D_{k+1}$ its closed neighborhood in $G_{k+1}$,
$d_{k+1}=|D_{k+1}|$, and $L_{i_{k+1}}$ one of the hyperedges containing  $u_{k+1}$, different from  $L_{i_1},\dots, L_{i_k}$.

We claim that this algorithm only stops after $a$ steps thus supplying the desired independent set $\{u_1,\dots,u_a\}$.
If we cannot define $u_{k+1}$ for a $k<a$, then $V_{k+1}$ is empty, and $|D_1|+\dots +|D_k|\geq s$.
As $D_1,\dots,D_k$ are non-empty, disjoint sets we get the contradiction
\[ 2|E(G)|\geq \sum_{u\in D_1\cup\dots\cup D_k} \deg(u)
        \geq d_1(d_1-1)+ \dots + d_k(d_k-1)
        \geq 2t(s,k) \geq 2 t(s,a-1). \quad {\Box}
        \]

%%% \vspace{3mm}
\bigskip
{\em Proof of Theorem~\ref{split}.}
Again, we will prove the upper bound (with $c_2=3$) by presenting a sufficient %c
 function $f$.
Let
\[f(v)=\left\{ \begin{array}{ll}
         s & \mbox{for $v\in A$},\\
         2 & \mbox{for $v\in Q$},\end{array} \right. \]
where $s:= \lfloor 3 \sqrt{(a-1)q} \rfloor $.
Note that $s\geq a$ and $t(s,a-1)>q$.
Consider any $f$-assignment, it provides sets $L_1, \dots, L_a$ and a graph $G$ with $|E(G)|=q$.
This assignment is sufficient %c
 by Theorem~\ref{turan}.

For the lower bound, we proceed as in Section~\ref{lower} and first consider
 $f$'s with $f(v)=2$ for all $v\in Q$.
Let $\{v_1,\dots,v_a\}$ be the vertices of $A$ and let
\[f(v)=\left\{ \begin{array}{ll}
         s_i & \mbox{for $v=v_i\in A$ and}\\
         2 & \mbox{for $v\in Q$} \end{array} \right. \]
such that $s_1\leq \dots \leq s_a$.
We claim that if $f$ is sufficient, %c
 then $q< t(s_i,i-1)$ for all $i\geq 2$.
These inequalities imply that $s_i\geq \sqrt{2(i-1)q}$ so
$$
  \sum f(v)\geq 2q+ \left( \sum_{1\leq i\leq a} \sqrt{i-1}\right)\sqrt{2q} \geq 2q+ \frac{1}{2} a\sqrt{(a-1)q}. $$
Suppose, on the contrary, that $q\geq t(s_i,i-1)$ for some $i\geq 2$.
Let $G$ be the graph consisting of $i-1$ disjoint cliques with sizes as equal as possible.
Assign the pairs corresponding to the edges of $G$ as the lists for the vertices of $Q$ and let $L(v_1)\subset \dots \subset L(v_i)=V(G)$.
This assignment is not sufficient. %c

To finish the proof of the lower bound for an arbitrary sufficient %c
 $f$ with
 $\sum f(v)=2q+as$ we use the inequality
(\ref{eq:q1+q3}) to obtain that $f(v)=2$ for all but at most $q-as$ vertices $v\in Q$.
Then we conclude the proof with an argument analogous to (\ref{tooMany}).
 %%% in the proof of Theorem~\ref{lb}.
 The details are omitted.  %%% conclude the proof for general lists.
\hfill{$\Box$}

\begin{prob} Suppose that
\[f(v)=\left\{ \begin{array}{ll}
         s_i & \mbox{for $v=v_i\in A$ and}\\
         2 & \mbox{for $v\in Q$} \end{array} \right. \]
such that $s_1\leq \dots \leq s_a$.

What conditions are sufficient and necessary for $f$ being $G_{a,q}$-sufficient? %c
\end{prob}

We already have seen that $q< t(s_i,i-1)$ for all $i\geq 2$ are necessary.
It is easy to see that $q\geq {{s_1}\choose{2}} +s_1(s_2-s_1)$ is also necessary,
 let $G$ be the graph that we get by taking a clique of order $s_2$
 and deleting edges of a clique of order $s_2-s_1$.
Assign its edges as the lists of $Q$, and let $L(v_1)\subset L(v_2)=V(G)$.
One is tempted to conjecture that these conditions altogether are already sufficient.

This is the same problem as to ask that how much the sizes of $L_i$'s can be decreased in Theorem~\ref{turan}.

\end{document}